\documentclass[11pt]{amsart}
\usepackage{amsmath}
\usepackage[all]{xy}
\usepackage{amssymb}
\usepackage{mathrsfs}
\usepackage{cite}
\usepackage{hyperref}
\usepackage{amsfonts}
\usepackage{mathdots}

\theoremstyle{plain}
\newtheorem{thm}{Theorem}[section]
\newtheorem{prop}[thm]{Proposition}
\newtheorem{lem}[thm]{Lemma}
\newtheorem{cor}[thm]{Corollary}

\theoremstyle{remark}
\newtheorem{rem}[thm]{Remark}

\theoremstyle{definition}
\newtheorem{defn}[thm]{Definition}

\newcommand{\cA}{\mathcal{A}}
\theoremstyle{conjecture}

\newcommand{\Ku}{\mathcal{K}u}

\def\HH{\mathrm{HH}}
\def\Hom{\mathrm{Hom}}
\def\Id{\mathrm{Id}}
\def\Jac{\mathrm{Jac}}
\DeclareMathOperator\oh{\mathcal{O}}

\title{Serre algebra, matrix factorization and categorical Torelli theorem for hypersurfaces}

\date{}

\begin{document}

\author{Xun Lin and Shizhuo Zhang}

\begin{abstract}
 Let $X$ be a smooth Fano variety. We attach a bi-graded associative algebra $\mathrm{HS}(\Ku(X))=\bigoplus_{i,j\in \mathbb{Z}} \Hom(\Id,S_{\Ku(X)}^{i}[j])$ to the Kuznetsov component $\Ku(X)$ whenever it is defined. Then we construct a natural sub-algebra of $\mathrm{HS}(\Ku(X))$ when $X$ is a Fano hypersurface and establish its relation with Jacobian ring $\mathrm{Jac}(X)$. As an application, we prove a categorical Torelli theorem for Fano hypersurface $X\subset\mathbb{P}^n(n\geq 2)$ of degree $d$ if $\mathrm{gcd}((n+1),d)=1.$ In addition, we give a new proof of the main theorem \cite[Theorem 1.2]{pirozhkov2022categorical} using a similar idea. 
\end{abstract}


\address{Max Planck Institute for Mathematics, Vivatsgasse 7, 53111 Bonn, Germany}
\email{xlin@mpim-bonn.mpg.de, lin-x18@tsinghua.org.cn}

\address{Max Planck Institute for Mathematics, Vivatsgasse 7, 53111 Bonn, Germany}
\address{Institut de Mathématiqes de Toulouse, UMR 5219, Université de Toulouse, Université Paul Sabatier, 118 route de
Narbonne, 31062 Toulouse Cedex 9, France}
\email{shizhuozhang@mpim-bonn.mpg.de,shizhuo.zhang@math.univ-toulouse.fr}

\maketitle

\section{Introduction}
Let $X$ be a smooth complex projective variety. Reconstruction of $X$ from its categorical invariant originates from Gabriel's thesis \cite{gabriel1962categories}, where the author proves $X$ can be recovered from its category of coherent sheaves. Later on, this theorem is generalized to arbitrary quasi-separated scheme in \cite{rosenberg1998noncommutative}. In the celebrated work \cite{bondal2001reconstruction}, the authors prove smooth Fano variety $X$ can be reconstructed from its bounded derived category $D^b(X)$ of coherent sheaves. 
Since the last decades, people are interested in reconstruction of Fano varieties from the non-trivial semi-orthogonal component $\Ku(X)$, known as Kuznetsov component, of its bounded derived category of coherent sheaves, called \emph{Categorical Torelli problem}. The first result in this direction is given in \cite{bernardara2012categorical}, where the authors prove categorical Torelli theorem for smooth cubic threefolds. Since then tremendous work have been carried out along this direction, see \cite{pertusi2023categorical} for a review of known results. In this paper, we focus on Fano hypersurfaces $X\subset\mathbb{P}^n$ of degree $d\leq n$. Our work is inspired by the paper \cite{huybrechts2016hochschild}, where the authors relate a variant of Hochschild cohomology $\mathrm{HH}(\Ku(X),(1))$ of Kuznetsov component $\Ku(X)$ to the Jacobian ring $J(X)$ of the hypersurface $X$, in particular they show the Hochschild cohomology ring is isomorphic to the Jacobian ring whenever $\Ku(X)$ is a Calabi-Yau category, hence establish a categorical Torelli theorem for cubic fourfolds. In addition, they suggested using the category of graded matrix factorizations $\operatorname{Inj_{coh}}(\mathbb{A}^{n+1},\mathbb{C}^{\ast},\omega)$ of a Fano hypersurface defined by a polynomial $\omega$ to reconstruct Jacobian ring via Hochschild cohomology $\mathrm{HH}^{\bullet}(\operatorname{Inj_{coh}}(\mathbb{A}^{n+1},\mathbb{C}^{\ast},\omega))$, which motivates our approach to this problem. On the other hand, for any smooth DG category $\cA$, there is a natural associative algebra $\mathrm{HS}(\cA))$ attached to it, called \emph{Serre algebra}(cf. Definition~\ref{def_Serre_algebra}), which is a Morita invariant of $\cA$. In \cite{bondal2001reconstruction}, the authors construct a subring of $\cA$, i.e. canonical ring to reconstruct smooth complex projective variety with canonical bundle ample or anti-ample. It is interesting to ask if certain sub-algebra of the Serre algebra $\mathrm{HS}(\Ku(X))$ of a smooth Fano variety $X$ can be used to determine the isomorphism class. In this article, we hope to answer this question. 

\subsection{Main results}
Let $X,X'\subset\mathbb{P}^n(n\geq 2)$ be Fano hypersurfaces of degree $d\leq n$. Instead of making additional assumption that the equivalence $\Ku(X)\simeq\Ku(X')$ is compatible with degree shifting functor $(1)$ and passing the equivalence $\operatorname{Inj_{coh}}(\mathbb{A}^{n+1},\mathbb{C}^{\ast},\omega)\simeq \operatorname{Inj_{coh}}(\mathbb{A}^{n+1},\mathbb{C}^{\ast},\omega')$ to the equivalence $\operatorname{Inj_{coh}}(\mathbb{A}^{n+1},\mathbb{C}^{\ast},\omega)/(1)\simeq \operatorname{Inj_{coh}}(\mathbb{A}^{n+1},\mathbb{C}^{\ast},\omega')/(1)$, we only assume that there is an equivalence $\Phi:\Ku(X)\simeq\Ku(X')$ and note that it commutes with Serre functors of $\Ku(X),\Ku(X')$ respectively. Then it is not hard to show that the associated \emph{Serre algebra}(cf. Definition~\ref{def_Serre_algebra}) of $\Ku(X)$ and $\Ku(X')$ are isomorphic. We construct a sub-algebra of Serre algebra $\mathrm{HS}(\Ku(X))$ and establish its relation with the Jacobian ring. 

\begin{thm}\label{main_theorem_sub_Serre_algebra}
Let $X\subset\mathbb{P}^n$ be a smooth hypersurface of degree $d\leq n$ defined by $\omega$. Consider its affine LG model whose underlying stack is $\mathcal{Z}=[\mathbb{A}^{n+1}/\mathbb{Z}_d]$ and associate $dg$ category of matrix factorization  $\operatorname{Inj_{coh}}(\mathbb{A}^{n+1},\mathbb{C}^{\ast},\omega)\simeq\Ku(X)$. Assume $\operatorname{gcd}((n+1),d)=1$. There is a $\mathbb{Z}$-graded sub-algebra $\bigoplus_{t\geq 0}\Hom(\Delta,\Delta(t))$ of the Serre algebra $\mathrm{HS}(\operatorname{Inj_{coh}}(\mathbb{A}^{n+1},\mathbb{C}^{\ast},\omega))$ such that 
$$\mathrm{Hom}(\Delta,\Delta(t))=\begin{cases}
    \mathrm{Jac}(\omega)_t, & n=2m,m\in\mathbb{Z}\\
    \mathrm{Jac}(\omega)_t, & n=2m+1,m\in\mathbb{Z}, t\neq \frac{(d-2)(n+1)}{2}\\
    \mathrm{Jac}(\omega)_{t}\oplus k^{d-1}, & n=2m+1, m\in\mathbb{Z},t=\frac{(d-2)(n+1)}{2}.
     \end{cases}
    $$
In particular, if $n$ is odd, the Serre algebra has a sub-graded algebra $\mathrm{Jac}(\omega)$ and if $n$ is even, the Serre algebra has a degree $d$ graded subspace $\mathrm{Jac}(\omega)_d$. 

%\end{enumerate}

\end{thm}

\begin{rem}\label{rem_grading_intro}
Since $\mathrm{gcd}((n+1),d)=1$, then there exists $k_1,k_2\in\mathbb{Z}$ such that $k_1(n+1)+dk_2=1$. Then it is shown in Proposition~\ref{subalgebra} that $S^{-k_1t}[2k_2t+(n+1)k_1t]\cong\Delta(t)$, for all $t\in\mathbb{Z}$, where $S$ is the Serre functor of $\operatorname{Inj_{coh}}(\mathbb{A}^{n+1},\mathbb{C}^{\ast},\omega)$ and $\Delta(1)$ is the twisting functor. Thus the $\mathbb{Z}$-graded piece $\mathrm{Hom}(\Delta,\Delta(t))\cong\mathrm{HS}(\operatorname{Inj_{coh}}(\mathbb{A}^{n+1},\mathbb{C}^{\ast},\omega))(-k_1t,2k_2t+k_1(n+1)t)$ and $\bigoplus_{t\geq 0}\mathrm{Hom}(\Delta,\Delta(t))$ is exactly the $\mathbb{Z}_{\geq 0}$-graded sub-algebra $\bigoplus_{t\in\mathbb{Z}_{\geq 0}}\mathrm{HS}(\operatorname{Inj_{coh}}(\mathbb{A}^{n+1},\mathbb{C}^{\ast},\omega))(-k_1t,2k_2t+k_1(n+1)t)$. 
\end{rem}

\begin{rem}
  Consider a smooth hypersuface of degree $d\leq \sum^{n+1}_{j=1}q_{j}-1$ in weighted projective space $\mathbb{P}(q_{1}, q_{2}\cdots, q_{n+1})$ with $\operatorname{gcd}(q_{1}, q_{2}, \cdots, q_{n+1})=1$. One can prove the same statement. Namely if $\operatorname{gcd}(\sum^{n+1}_{j=1}q_{j}, d)=1$, then there is a sub-algebra $\bigoplus_{t\geq 0}\Hom(\Delta,\Delta(t))$ of $\mathrm{HS}(\operatorname{Inj_{coh}}(\mathbb{A}^{n+1},\mathbb{C}^{\ast},\omega)),$ where $\mathbb{C}^{\ast}$-action on $\mathbb{A}^{n+1}$ is of weight $(q_{1}, q_{2}, \cdots, q_{n+1})$. 
  %such that   
%$$\mathrm{Hom}(\Delta,\Delta(t))=
 %   \mathrm{Jac}(\omega)_t, n=2m, m\in\mathbb{Z}
 %   $$
\end{rem}
As an application, we establish \emph{Categorical Torelli theorem} for smooth hypersurface $X\subset\mathbb{P}^n$ of degree $d$ whenever $\mathrm{gcd}((n+1),d)=1$. Namely, we show the Kuznetsov component $\Ku(X)$(cf.Proposition~\ref{prop_Orlov_identification})
determines its isomorphism class. 

\begin{thm}\label{main_theorem_second}
   Let $X,X'\subset\mathbb{P}^n$ be degree $d$ smooth hypersurfaces. Assume $\mathrm{gcd}(n+1,d)=1$. If there is a Fourier-Mukai equivalence $\Ku(X)\simeq \Ku(X')$, then $X\cong X'$.   
\end{thm}
\begin{rem}
  Consider a smooth hypersufaces of degree $d\leq \sum^{n+1}_{j=1}q_{j}-1$ in weighted projective space $\mathbb{P}(q_{1}, q_{2}\cdots, q_{n+1})$ with $\operatorname{gcd}(q_{1}, q_{2}, \cdots, q_{n+1})=1$. By similar arguments in the proof of Theorem~\ref{main_theorem_second}(or Theorem~\ref{theorem_categorical_Torelli_thm_without_rotation_functors}), one is able to obtain \emph{Categorical Torelli theorem} for a series of hypersurfaces. In particular, our method works for a 
  degree $6$ hypersurface in weighted projective space $\mathbb{P}(1,1,1,1,3)$, which is isomorphic to an index one degree two Fano threefold(See Section~\ref{section_appendix}).\footnote{In an earlier version of the paper, a weaker version of Theorem~\ref{main_theorem_second} is obtained. Namely the statement is obtained if $\mathrm{gcd}(2(n+1),d)=1$ for non-weighted case and $\mathrm{gcd}(2\Sigma_{j=1}^{n+1}q_j,d)=1$ for weighted case, while we are informed by Paolo Stellari that in their paper \cite{Lahoz2023categorical}, categorical Torelli theorem for sextic hypersurface in $\mathbb{P}(1,1,1,1,3)$ is established. Then we found that our method applies to their case and a stronger statement (as presented in the paper) is obtained and we give a new proof of the sextic weighted hypersurface case in Appendix.}
\end{rem}
Let us briefly explain the idea of the proof. We work on $dg$ category of graded matrix factorization  $\operatorname{Inj_{coh}}(\mathbb{A}^{n+1},\mathbb{C}^{\ast},\omega)$. An equivalence $\Phi:\Ku(X)\simeq\Ku(X')$ in $\operatorname{Hqe(dg-cat)}$ induces an equivalence $\operatorname{Inj_{coh}}(\mathbb{A}^{n+1},\mathbb{C}^{\ast},\omega)\simeq \operatorname{Inj_{coh}}(\mathbb{A}^{n+1},\mathbb{C}^{\ast},\omega')$, commuting with Serre functors on both categories. Whenever $n$ is even, we get graded ring isomorphism $\mathrm{Jac}(\omega)\cong\mathrm{Jac}(\omega')$ coming from a degree one map, then Mather-Yau reconstruction theorem \cite[Proposition 1.1]{Don} gives $X\cong X'$. If $n$ is odd, we get the isomorphism of degree $d$ or $d-1$ components of $\mathrm{Jac}(\omega)$. Then similar arguments shows $X\cong X'$. 

As a corollary, we give some interesting examples where \emph{Categorical Torelli theorem} holds. 

%another proof of \emph{Categorical Torelli theorem} for smooth cubic threefolds established in \cite{bernardara2012categorical} and produce a new case: cubic surfaces. 

\begin{cor}\label{main_theorem_third}
Categorical Torelli theorem holds for following Fano varieties:
\leavevmode\begin{enumerate}
    \item Cubic hypersurfaces of dimension $3k-1$ and $3k$ with $k\geq 1$. 
    \item Quintic fourfolds.
\end{enumerate}
\end{cor}

In the paper \cite{pirozhkov2022categorical}, the author shows in \cite[Theorem 1.2]{pirozhkov2022categorical} that a class of Fano hypersurfaces $X$ are determined by the Kuznetsov components $\Ku(X)$ together with rotation functors $(1):\Ku(X)\simeq\Ku(X)$, 
which generalizes a result in \cite[Corollary 2.10]{huybrechts2016hochschild}. Using the framework of matrix factorization, we give a simple proof for \cite[Theorem 1.2]{pirozhkov2022categorical}. 

\begin{thm}\label{main_theorem_fourth}
  Let $X$ and $X'$ be smooth hypersurfaces of degree $d<n+1$ in $\mathbb{P}^{n}(n\geq 2)$.  If there is a Fourier-Mukai equivalence of pairs $(\Ku(X),(1))\cong (\Ku(X'),(1)')$ , then $X\cong X'$.   \end{thm}

\subsection{Related Work}
Our work is inspired by the paper \cite{huybrechts2016hochschild}, where the authors suggest using the framework of graded matrix factorization of a smooth hypersurface to relate the extended Hochschild cohomology $\mathrm{HH}(\Ku(X),(1))$ to the Jacobian ring $J(X)$, but they do all the work on derived category side. In \cite{pirozhkov2022categorical}, the author generalized results \cite[Corollary 2.10]{huybrechts2016hochschild} to arbitrary Fano hypersurfaces. In an upcoming paper \cite{Rennemo2023hochschild}, the author proves if $d$ does not divide $n+1$ and the pair $(d,n)$ is not of the form $(4,4k+1)$, then the Kuznetsov component alone reconstructs $X$. In upcoming work \cite{Lahoz2023categorical} and \cite{Dell2023categorical}, the authors prove categorical Torelli theorem for sextic hypersurface in weighted projective space $\mathbb{P}(1,1,1,1,3)$ via completely different methods. 

\begin{rem}
While preparing the paper, we learned that the authors of the paper  \cite{belmans2023hochschild} also define and study the Serre algebra under the name Hochschild-Serre cohomology in \cite[Definition A.1]{belmans2023hochschild}, where they give a formula for Hochschild cohomology of Hilbert scheme of points on a surface in terms of Hochschild-Serre cohomology(as a bi-graded vector space)on the surface. 
\end{rem}

%observes the action of $\mathrm{HH}^*(\cA_X)$ on $\mathrm{HH}_*(\cA_X)$ recovers some graded pieces of the Jacobian ring, enough to recover $X$ via Donagi's proof of the generic hypersurface Torelli theorem. 

%In \cite{orlov2003derived}, the author defined the Hochschild algebra $\mathrm{HA}$ by letting $\Ku(X)$ be the whole bounded derived category $D^b(X)$.

\subsection{Organization of the article}
In Section~\ref{section_intro_matrix_factorizations} we introduce the terminology of graded matrix factorization associated with a hypersurface in projective space. Then we describe an important auto-equivalence on the category of matrix factorization $\operatorname{Inj_{coh}}(\mathbb{A}^{n+1},\mathbb{C}^*,\omega)$. In Section~\ref{section_Serre_algebra}, we introduce an associative algebra naturally attached to any smooth and proper differential graded($dg$) category, called \emph{Serre algebra} and show it is a Morita invariant. Then we give several examples of Serre algebra for various $dg$ category, in particular, we construct interesting sub-algebra of Serre algebra for category of matrix factorization $\operatorname{Inj_{coh}}(\mathbb{A}^{n+1},\mathbb{C}^*,\omega)$ corresponding to a Fano hypersurface in $\mathbb{P}^n$, proving Theorem~\ref{main_theorem_sub_Serre_algebra}. In Section~\ref{theorem_categorical_Torelli_thm_without_rotation_functors} we prove Theorem~\ref{main_theorem_second} and Corollary~\ref{main_theorem_third}. In Section~\ref{categorical_torelli_with_rotation_functors} we prove Theorem~\ref{main_theorem_fourth}. In Section~\ref{section_appendix}, we give a new proof of \emph{Categorical Torelli theorem} for sextic hypersurface in $\mathbb{P}(1,1,1,1,3)$ as recently shown in \cite{Lahoz2023categorical}. 

\subsection{Acknowledgement}
We would like to thank Will Donovan, Daniel Huybrechts, Ziqi Liu, J{\o}rgen Rennemo, Ed Segal, Junwu Tu, and Jieheng Zeng for useful conversations on related topics. We thank Paolo Stellari for letting us know about the paper \cite{Lahoz2023categorical}, where they prove categorical Torelli theorem for degree $6$ hypersurfaces in weighted projective space $\mathbb{P}(1,1,1,1,3)$
via a completely different method and the same statement is also proved in another upcoming preprint \cite{Dell2023categorical}. We also would like to thank the referee for the careful reading and for providing detailed comments. SZ is supported by ANR project FanoHK, grant ANR-20-CE40-0023, Deutsche Forschungsgemeinschaft under Germany's Excellence Strategy-EXC-2047$/$1-390685813, and partially supported by GSSCU2021092. Part of the work was finished when XL and SZ are visiting Max-Planck and Hausdorff institutes for Mathematics. 
They are grateful for excellent working conditions and hospitality.

\subsection{Data availability}
Data sharing does not apply to this article as no datasets were generated or analyzed during the current research.

\subsection{Conflict of interest}
On behalf of all authors, the corresponding author states that there is no Conflict of interest.

\section{dg category of graded matrix factorizations}\label{section_intro_matrix_factorizations}
In this section, we recall the terminology of $dg$-category of matrix factorization. We follow the context in \cite{BFK}. We refer the reader to \cite{keller2006differential} for the basic $dg$ categories. 
Let $\operatorname{Hqe(dg-cat)}$ be the localizing of $\operatorname{dg-cat}$ with respect to the quasi-equivalences of $dg$ categories.   Let $(X,G,L,\omega)$ be a 
quadruple where $X$ is a quasi-projective variety with $G$ action, $G$ is a reductive algebraic group, $L$ is a $G$-equivariant line bundle and $\omega$ is a $G$-invariant section of $L$. Our main example is $(\mathbb{A}^{n+1},\mathbb{C}^{\ast},\mathcal{O}(d),\omega)$. The action of $\lambda\in\mathbb{C}^{\ast}$ on $\mathbb{A}^{n+1}$ is given by $\lambda\cdot(x_{0},x_{1},\cdots, x_{n})=(\lambda\cdot x_{0},\lambda\cdot x_{1},\cdots, \lambda\cdot x_{n})$. $\omega$ is a $\mathbb{C}^{\ast}$-invariant section of $\mathcal{O}(d)$. Namely $\omega$ is a degree $d$ polynomial. We always assume $\omega$ has only isolated singularity at $0\in \mathbb{A}^{n+1}$. 

We have $dg$ category $\operatorname{Fact}(X,G,L,\omega)$, whose objects are a quadruple $(\mathcal{E}_{-1},\mathcal{E}_{0},\Phi_{-1},\Phi_{0})$, where $\mathcal{E}_{-1}$ and $\mathcal{E}_{0}$ are $G$-equivariant quasi-coherent sheaves, $\Phi_{-1}:\mathcal{E}_{0}\rightarrow \mathcal{E}_{-1}\otimes L$ and $\Phi_{0}:\mathcal{E}_{-1}\rightarrow \mathcal{E}_{0}$ are morphism of $G$-equivariant sheaves such that
\begin{align*}
    \Phi_{-1}\circ \Phi_{0}=\omega.\\
    (\Phi_{0}\otimes L)\circ \Phi_{-1}=\omega.
\end{align*}
The space of morphisms in $\operatorname{Fact}(X,G,L,\omega)$ are the internal Hom of $G$-equivariant sheaves while extending the pairs of morphisms to certain $\mathbb{Z}$-graded complexes. We point out the reference \cite{BFK} for interested reader. Let $\operatorname{Inj}(X,G,L,\omega)\subset \operatorname{Fact}(X,G,L,\omega)$ be a dg sub-category whose components are $G$-equivariant injective quasi-coherent sheaves. There is a category $\operatorname{Acycli}(\operatorname{Fact}(X,G,L,\omega))$ which imitates acyclic complexes in category of complexes of sheaves. The absolute derived category $D^{abs}(\operatorname{Fact}(X,G,L,\omega))$ is the homotopy category of $dg$ quotient $\frac{\operatorname{Fact}(X,G,L,\omega)}{\operatorname{Acyclic}(\operatorname{Fact}(X,G,L,\omega))}$ in $\operatorname{Hqe(dg-cat)}$. 
\begin{lem}
   The natural morphism $\operatorname{Inj}(X,G,L,\omega)\rightarrow D^{abs}(\operatorname{Fact}(X,G,L,\omega))$
   induces isomorphism in homotopic categories.
\end{lem}
Let $\operatorname{Inj_{coh}}(X,G,L,\omega)\subset \operatorname{Inj}(X,G,L,\omega)$ be a $dg$ sub-category whose objects are quasi-isomorphic to objects with coherent components in category $\operatorname{Fact}(X,G,L,\omega)$.
%\begin{defn}
  
  Define shifting functor $$[1]: (\mathcal{E}_{-1},\mathcal{E}_{0},\Phi_{-1},\Phi_{0})\mapsto (\mathcal{E}_{0},\mathcal{E}_{-1}\otimes L,-\Phi_{0},-\Phi_{-1}\otimes L).$$ 
%\end{defn}
With cone construction, the homotopic categories $[\operatorname{Inj_{coh}}(X,G,L,\omega)]$ is a triangulated category which is isomorphic to graded matrix factorization in \cite{orlov2009derived} for $(\mathbb{A}^{n+1},\mathbb{C}^{\ast},\mathcal{O}(d),\omega)$.
%\begin{defn}
  
Denote by  
  $$\{1\}:\operatorname{Inj_{coh}(\mathbb{A}^{n+1},\mathbb{C}^{\ast},\mathcal{O}(d),\omega)}\rightarrow \operatorname{Inj_{coh}(\mathbb{A}^{n+1},\mathbb{C}^{\ast},\mathcal{O}(d),\omega)}$$
  the twisting functor which maps 
$$\xymatrix{\mathcal{E}_{-1}
\ar[r]^{\Phi_{0}}&\mathcal{E}_{0}\ar[r]^{\Phi_{-1}}&\mathcal{E}_{-1}(d)}$$
to 
$$\xymatrix{\mathcal{E}_{-1}(1)
\ar[r]^{\Phi_{0}(1)}&\mathcal{E}_{0}(1)\ar[r]^{\Phi_{-1}(1)}&\mathcal{E}_{-1}(d+1)}$$
%\end{defn}
%The following theorem is established in \cite[Theorem 2.5]{orlov2009derived}. 
Clearly, we have equality of functors $\{d\}:=\{1\}^{d}=[2]$.
\par
 Let $X\subset \mathbb{P}^{n}$ be a smooth hypersurface of degree $d\leq n$ defined by $\omega$. Let 
  $$\Ku(X):=\Big\langle \mathcal{O}_{X},\mathcal{O}_{X}(1),\cdots, \mathcal{O}_{X}(n-d)\Big\rangle^{\perp}.$$
 Consider the natural enhancement $\operatorname{Inj_{coh}}(X)$, and let $\Ku_{dg}(X)$ be a $dg$ subcategory that enhance $\Ku(X)$. Write $(1)$ as the quasi-endofunctor (Fourier-Mukai type) of $\Ku_{dg}(X)$ that define degree shifting functor $(1):\Ku(X)\rightarrow \Ku(X)$ in the sense of \cite{huybrechts2016hochschild}.
\begin{thm}\label{prop_Orlov_identification}\cite[Theorem 6.13]{BFK}
  There is an equivalence in $\operatorname{Hqe(dg-cat)}$,
  $$\Phi:\operatorname{Inj_{coh}}(\mathbb{A}^{n+1},\mathbb{C}^{\ast},\omega)\cong \Ku_{dg}(X).$$
  In particular, there is an isomorphism of quasi-funtors
  $$\Phi\circ\{1\}\cong (1)\circ \Phi,$$
  where $\Phi$ is a quasi-functor $\operatorname{Inj_{coh}}(\mathbb{A}^{n+1},\mathbb{C}^{\ast},\mathcal{O}(d),\omega)\rightarrow\mathrm{Inj}_{\mathrm{coh}}(X)$, which is defined explicitly in the proof of \cite[Theorem 6.13]{BFK}.
\end{thm}
\begin{proof}
 Firstly, we have quasi-functor
 $$\Phi: \operatorname{Inj_{coh}}(\mathbb{A}^{n+1},\mathbb{C}^{\ast},\mathcal{O}(d),\omega)\rightarrow \Ku_{dg}(X),$$ 
 and quasi-functor 
 $$\Phi^{!}: \operatorname{Inj_{coh}}(X)\rightarrow \operatorname{Inj_{coh}}(\mathbb{A}^{n+1},\mathbb{C}^{\ast},\mathcal{O}(d),\omega)$$ 
 such that $\Phi^{!}\circ \Phi\simeq\Id$. Since $\Phi$ defines an equivalence of triangulated category $[\operatorname{Inj_{coh}}(\mathbb{A}^{n+1},\mathbb{C}^{\ast},\mathcal{O}(d),\omega)]\simeq \Ku(X)$, and both $\operatorname{Inj_{coh}}(\mathbb{A}^{n+1},\mathbb{C}^{\ast},\mathcal{O}(d),\omega)$ and $\Ku_{dg}(X)$ are exact $dg$ categories, therefore $\Phi$ is an isomorphism in $\operatorname{Hqe(dg-cat)}$.
 The equality $\Phi^{!}\circ \Phi=\Id$ implies $\Phi^{!}$ is the inverse of $\Phi$ when restricting to $\Ku_{dg}(X)$. Write $(1)_{F}:\operatorname{Inj_{coh}}(X)\rightarrow \operatorname{Inj_{coh}}(X)$ as the Fourier-Mukai functor that restricts to rotation functor $(1)$ on $\Ku_{dg}(X)$ in \cite[Theorem 6.13]{BFK}. The isomorphims of quasi-functors $\Phi^{!}\circ(1)_{F}\circ \Phi\simeq \{1\}$ implies an isomorphism of quasi-functors $\Phi^{-1}\circ (1)\circ \Phi\simeq \{1\}$.
\end{proof}

\section{Serre algebra}\label{section_Serre_algebra}

\begin{thm}\label{functorbimodule}\cite{toen2007homotopy} Let $\mathcal{A}$ be a $dg$ category over the field $k$. In $\operatorname{Hqe(dg-cat)}$, we have isomorphism, 
$$\mathcal{D}_{dg}(\mathcal{A}^{op}\otimes \mathcal{A})\cong \mathcal{R}\Hom_{c}(\mathcal{D}_{dg}(\mathcal{A}),\mathcal{D}_{dg}(\mathcal{A})),$$
where $\mathcal{R}\Hom_{c}$ is the quasi-functor preserving coproduct.
\end{thm}
If $\mathcal{A}$ is a smooth proper $dg$ category, then the bimodules for Serre functor and inverse of Serre functor are constructed explicitly in \cite{shklyarov2007serre}, where the author defines Serre functor(inverse Serre functor) of the triangulated category $\operatorname{Perf}(\mathcal{A})$ in the usual sense. From now on, we interpret 
those bimodules by quasi-functors by Theorem~\ref{functorbimodule}. 

%$$HH^{m}(\mathcal{A})=Hom_{D(\mathcal{A}\otimes \mathcal{A}^{op})}(\mathcal{A},\mathcal{A}[m]).$$
%$$HH_{m}(\mathcal{A})=Hom_{D(\mathcal{A}\otimes \mathcal{A}^{op})}(\mathcal{A},\mathcal{A}^{\ast}[m]).$$
\begin{defn} The Hochschild (co)homology of a smooth proper $dg$ category $\mathcal{A}$ are defined as,
    $$\HH^{m}(\mathcal{A})=\Hom(\Id,\Id[m]),$$
$$\HH_{m}(\mathcal{A})=\Hom(\Id,S[m]),$$
\end{defn}
where $S$ is the Serre functor of $\cA$. The Hochschild cohomology is an algebra, and the homology is a graded module over the Hochschild homology. We define an algebra that contains Hochschild cohomology and Hochschild homoloy and encodes the algebra structure of Hochschild cohomology and the module structure of Hochschild homology over the Hochschild cohomology.

\begin{defn}(Serre algebra)\label{def_Serre_algebra}
   Let $\mathcal{A}$ be a smooth proper $dg$ category and $S$ the Serre functor of $\cA$.  Define bi-graded algebra 
   $$\mathrm{HS}(\cA):=\bigoplus_{m,n\in\mathbb{Z}}\mathrm{HS}(\cA)(m,n)=\bigoplus_{m,n\in \mathbb{Z}} \Hom(\Id,S^{m}[n]).$$
 The multiplication map 
 $$\xymatrix@C=2cm{\Hom(\Id,S^{m_1}[n_1])\times \Hom(\Id,S^{m_2}[n_2])\ar[r]^{\times}&\Hom(\Id,S^{m_1+m_2}[n_1+n_2])}$$
 is defined as follows.
 For element $(a,b)\in \Hom(\Id,S^{m_1}[n_1])\times \Hom(\Id,S^{m_2}[n_2])$, $a\times b$ is defined as the composition
 $$\xymatrix@C=3.5cm{\Id\ar[r]^{b}&\Id\circ S^{m_2}[n_2]\ar[r]^{a\circ \Id}&S^{m_1}[n_1]\circ S^{m_2}[n_2]=S^{m_1+m_2}[n_1+n_2]
 }.$$
 We check the associativity. Namely for elements 
 $a\in \Hom(\Id,S^{m_1}[n_1])$, $b\in \Hom(\Id,S^{m_2}[n_2])$, and $c\in \Hom(\Id,S^{m_3}[n_3])$, we have 
 $(ab)c=a(bc)=abc$. For example, $a(bc)=abc$ follows from the commutative diagram,
 $$\xymatrix@C=2cm{\Id\ar[r]^{c}\ar[drr]_{bc}&S^{m_3}[n_3]\ar[r]^{b} &S^{m_3}[n_3]\circ S^{m_2}[n_2] \ar[r]^{a}\ar@{=}[d]&S^{m_3}[m_3]\circ S^{m_2}[n_2]\circ S^{m_1}[n_1]\ar@{=}[d]\\
 & &S^{m_3+m_2}[n_3+n_2]\ar[r]^{a} &S^{m_3+m_2+m_1}[n_3+n_2+n_1]}$$
\end{defn}
Let $\operatorname{Hmo(dg-cat)}$ be the localization of $\operatorname{dg-cat}$ with respect to the Morita equivalences of $dg$ categories. If $\mathcal{A}$ and $\mathcal{B}$ are smooth and proper, 
$\Hom_{\operatorname{Hmo(dg-cat)}}(\mathcal{A},\mathcal{B})$ is the isomorphism classes of $\operatorname{Perf}(\mathcal{A}^{op}\otimes \mathcal{B})$\cite[Corollary 1.44]{tabuada2015noncommutative}, and the composition corresponds to tensor product.

%\begin{lem}
%Let $\mathcal{A}$ and $\mathcal{B}$ be two smooth proper $dg$ categories, and $\Phi:\mathcal{A}\cong \mathcal{B}$ be the isomorphism in $Hmo(dg-cat)$, then $\Phi\circ S_{\mathcal{A}}\cong S_{\mathcal{B}}\circ \Phi$.
%\end{lem}
%\begin{proof}
%  We write $Z$ as the bimodule coresponds to $\Phi$.   
%\end{proof}
\begin{thm}\label{Moritaserrealgebra}
If $\mathcal{A}\cong \mathcal{B}$ in $\operatorname{Hmo(dg-cat)}$, then $\mathrm{HS}(\cA)\cong \mathrm{HS}(\mathcal{B})$.    
\end{thm}
\begin{proof}
  There is an isomorphism of Serre functor $S_{\mathcal{A}}\circ \Phi \cong \Phi\circ S_{\mathcal{B}}$. Hence, $\Phi$ induces isomorphism for any integer $m$, $n$,
  $$\Hom(\Id_{\mathcal{A}},S_{\mathcal{A}}^{m}[n])\cong \Hom(\Id_{\mathcal{B}},S_{\mathcal{B}}^{m}[n]).$$
  by the following commutative diagram,
  $$\xymatrix{\Id_{\mathcal{B}}\ar[r]\ar[d]^{\cong}&S^{m}_{\mathcal{B}}[n]\ar[d]^{\cong}\\
  \Phi^{-1}\circ \Id_{\mathcal{A}}\circ \Phi\ar[r]&\Phi^{-1}\circ S^{m}_{\mathcal{A}}[n]\circ\Phi}$$
  The isomorphism is an isomorphism of algebra since both algebra are defined by composition of functors.
\end{proof}

\subsection{Examples of Serre Algebras}
In this section, we give examples of Serre algebra for various categories. 

\subsubsection{Orlov's algebra $\mathrm{HA}(X)$}\label{section_orlov_algebra}
Let $\cA=D^b(X)$ be the bounded derived category of coherent sheaves on a smooth projective variety $X$. In this case, $S_{\cA}=-\otimes\omega_X[l]$, where $l=\mathrm{dim}X$. Thus the Serre algebra $\mathrm{HS}(\cA))$ is given by 
\begin{align*}
\mathrm{HS}(\cA):=\bigoplus_{m,n\in\mathbb{Z}}\mathrm{Hom}(\mathrm{Id},\mathrm{S}_{D^b(X)}^m[n])
\cong&\bigoplus_{m,n\in\mathbb{Z}}\mathrm{Hom}_{D^b(X\times X)}(\iota_*\oh_X,\iota_*\omega_X^{\otimes m}[ml+n])\\
\cong&\bigoplus_{m,n\in\mathbb{Z}}\mathrm{Ext}_{X\times X}^{ml+n}(\iota_*\oh_X,\iota_*\omega_X^{\otimes m}),
\end{align*}
where $\iota:X\hookrightarrow X\times X$ be the diagonal inclusion. It is clear that $\mathrm{HS}(\cA)$ is isomorphic to the bi-graded algebra $\mathrm{HA}(X)$ in \cite{orlov2003derived}. In particular, for a given $m \in \mathbb{Z}$, the graded piece $\mathrm{HS}(\cA)(m,-ml)$ equals
\[
\mathrm{Hom}_{X\times X}(\iota_*\oh_X,\iota_*\omega_X^{\otimes m}) = \mathrm{Hom}_X(\oh_X,\omega_X^{\otimes m}),
\]
and summing over $m \ge 0$ exhibits the canonical ring of $X$ as a subring of $\mathrm{HS}(\cA))$.

%In particular, if $ml+n=0$, then $\cA_S=\bigoplus_{m\geq 0}\mathrm{Hom}_{X\times X}(\iota_*\oh_X,\iota_*\omega_X^{\otimes m})\cong\bigoplus_{m\geq 0}\mathrm{Hom}_X(\oh_X,\omega_X^{\otimes m})$ is the canonical ring of $X$. 

%\subsubsection{Kuznetsov's component of quartic double solids}
%Let $Y$ be a smooth quartic double solid, whose Kuznetsov component $\Ku(Y)$ is an Enriques category. Its Serre functor $\mathrm{S}_{\Ku(Y)}\cong\tau\circ [2]$, where $\tau$ is the geometric involution induced by double cover. (To complete...)

\subsubsection{$dg$-category of matrix factorization on affine LG model}
According to Orlov' sigma/LG correspondence, the Kuznetsov components of hypersurfaces $X\subset\mathbb{P}^n$ of degree $d$ are affine LG model. First, the Serre functor of $\operatorname{Inj_{coh}}(\mathbb{A}^{n+1},\mathbb{C}^{\ast},\omega)$(Theorem~\ref{prop_Orlov_identification}) is $-\otimes \mathcal{O}_{\mathbb{A}^{n+1}}(-n-1)[n+1]$ \cite[Theorem 1.2]{FeveroKellyLGfractionalcy}. According to \cite{BFK}, the natural functors can be reinterpreted as Fourier-Mukai transformation of kernels, and the natural transformation between these functors is a morphism of kernels. We write $\Delta(m)$ as the kernel of functor $-\otimes\mathcal{O}_{\mathbb{A}^{n+1}}(m)$. Next, we recall a key theorem in \cite[Theorem 1.2]{BFK}. For $g\in\mu_{d}=\langle \exp{\frac{2\pi i}{d}}\rangle$, let $W_{g}$ be the conormal sheaf of $(\mathbb{A}^{n+1})^{g}$ in $\mathbb{A}^{n+1}$, note that $(\mathbb{A}^{n+1})^{g}$ is always a linear subspace of $\mathbb{A}^{n+1}$ and there we use the standard definition of Koszul cohomology. Let $k_{g}$ the character of $\operatorname{det}(W_{g})$. Let $N$ be the dimension of $\mathbb{A}^{g}$, and $\omega_{g}:=\omega\vert_{(\mathbb{A}^{n+1})^{g}}$.
Consider the Koszul complex $K_{\bullet}\langle \partial_{1}\omega_{g}, \partial_{2}\omega_{g}, \cdots, \partial_{N}\omega_{g}\rangle$, where $\langle \partial_{1}\omega_{g}, \partial_{2}\omega_{g}, \cdots, \partial_{N}\omega_{g}\rangle$ is a sequence in the coordinate ring $R$ of $\mathbb{A}^{g}$. It is the complex,
$$\xymatrix{\wedge^{N}R^{N}\ar[r]^{d}&\wedge^{N-1}R^{N}\ar[r]^{d}&\ar[r]\cdots\ar[r]&\wedge^{2}R^{N}\ar[r]^{d}&R^{N}\ar[r]^{d}&R}$$
with differential $d$ mapping the basis elements $e_{i_{1}}\wedge e_{i_{2}}\wedge\cdots \wedge e_{i_{p}}(1\leq i_{1}<i_{2}<\cdots<i_{p}\leq N)$ to 
$$d(e_{i_{1}}\wedge e_{i_{2}}\wedge\cdots \wedge e_{i_{p}})=\sum^{p}_{k=1}(-1)^{k+1}\partial_{i_{k}}\omega_{g}\cdot e_{i_{1}}\wedge \cdots \wedge \hat{e}_{i_{k}}\wedge \cdots \wedge e_{i_{p}}.$$
 We write $H^{\bullet}(d\omega_g)$ for the cohomology of the Koszul complex
 $K_{\bullet}\langle \partial_{1}\omega_{g}, \partial_{2}\omega_{g}, \cdots, \partial_{N}\omega_{g}\rangle$. Assume $\omega$ has a unique isolated singularity, then $H^{\bullet}(d\omega_{g})$ has only non-trivial cohomology at degree zero, and $H^{\bullet}(d\omega_{g})=H^{0}(d\omega_{g})=\mathrm{Jac}(\omega_{g})$.
 %\jorgen[inline]{Add a definition of ``Koszul cohomology''.}
\begin{thm}\cite[Theorem 5.9]{BFK}\label{extendedhochschild}
 Assume $\omega$ has an isolated singularity exactly at $0\in \mathbb{A}^{n+1}$, then 
 %\begin{align*}
 %\Hom(\Delta,\Delta(m)[t])\cong& (\bigoplus_{g\in \mu_{d},\ l\geq 0,\  t-\operatorname{rk}\omega_{g}=2u }H^{2l}(d\omega_{g})(m-k_{g}+d(u-l))\\
 %&\oplus \bigoplus_{g\in \mu_{d},\ l\geq 0,\  t-\operatorname{rk}W_{g}=2u+1 }H^{2l+1}(d\omega_{g})(m-k_{g}+d(u-l)))^{\mathbb{C}^{\ast}} 
 %\end{align*}
%Furthermore, since $H^{\bullet}(d\omega_{g})$ has only non-trivial cohomology at degree zero, namely $H^{\bullet}(d\omega_{g})=H^{0}(d\omega_{g})=\mathrm{Jac}(\omega_{g})$, we have
\begin{align*}
 \Hom(\Delta,\Delta(t)[m])\cong& (\bigoplus_{g\in \mu_{d},\  t-\operatorname{rk}W_{g}\  \textbf{is even} }\mathrm{Jac}(d\omega_{g})(t-k_{g}+d(\frac{m-\operatorname{rk}W_{g}}{2})))^{\mathbb{C}^{\ast}}. 
 \end{align*}
\end{thm}
\begin{rem}\label{compatibility}
%As $(\mathbb{A}^{n+1})^{g}=(0,\cdots,0)$ if $g\neq 1$, $Jac(\omega_{g})=k(0)$. In the cases $g\neq 1$ has no contribution to $\Hom(\Delta,\Delta(t_{1}))$ and $\Hom(\Delta,\Delta(t_{2}))$, 
If $\Hom(\Delta,\Delta(t_{1}))\cong \Jac(\omega)_{t_{1}}$ and $\Hom(\Delta,\Delta(t_{2}))\cong \Jac(\omega)_{t_{2}}$, then the multiplication 
$$\Hom(\Delta,\Delta(t_{1}))\times \Hom(\Delta,\Delta(t_{2}))\rightarrow \Hom(\Delta,\Delta(t_{1}+t_{2}))$$
is the composition of functions on $\mathbb{A}^{n+1}$ (namely the product of polynomials) while identifying with certain graded pieces of Jacobian algebra $\Jac(\omega)$.
\end{rem}

\begin{prop}\label{subalgebra}
Consider the affine $LG$ model $\operatorname{Inj_{coh}}(\mathbb{A}^{n+1},\mathbb{C}^{\ast},\omega)$. Assume $\operatorname{gcd}((n+1),d)=1$. There is a sub-algebra $\bigoplus_{t\geq 0}\Hom(\Delta,\Delta(t))$ of $\mathrm{HS}(\operatorname{Inj_{coh}}(\mathbb{A}^{n+1},\mathbb{C}^{\ast},\omega))$ such that 
$$\mathrm{Hom}(\Delta,\Delta(t))=\begin{cases}
    \mathrm{Jac}(\omega)_t, & n=2m,m\in\mathbb{Z}\\
    \mathrm{Jac}(\omega)_t, & n=2m+1,m\in\mathbb{Z}, t\neq \frac{(d-2)(n+1)}{2}.\\
    \Jac(\omega)_{t}\oplus k^{d-1}, & n=2m+1,m \in \mathbb{Z}, t=\frac{(d-2)(n+1)}{2}.
    \end{cases}
    $$
In particular, if $n$ is odd, the Serre algebra has a sub-graded algebra $\mathrm{Jac}(\omega)$ and if $n$ is even, the Serre algebra has a degree $d$ graded subspace $\mathrm{Jac}(\omega)_d$. 
\end{prop}
\begin{proof}
Firstly, It is known that the Serre functor $S\cong \Delta(-(n+1))[n+1]$ \cite[Theorem 1.2]{FeveroKellyLGfractionalcy}, and we have $\Delta(d)=[2]$. If $\operatorname{gcd}(n+1,d)=1$, then there exists $k_{1}, k_{2}\in\mathbb{Z}$ such that $k_{1}(n+1)+k_{2}d=1.$
Then 
$$S^{k_{1}}\cong\Delta(-k_{1}(n+1))[k_{1}(n+1)] \cong\Delta(-1+k_{2}d)[k_{1}(n+1)]\cong\Delta(-1)[2k_{2}+k_{1}(n+1)].$$  

Thus for any integer $t$, $(S^{-k_{1}}[2k_{2}+k_{1}(n+1)])^{t}\cong\Delta(t)$.
The composition of $-\otimes \mathcal{O}(-)$ is the same as the composition of Serre functors, hence $\bigoplus_{t\geq 1}\Hom(\Delta,\Delta(t))$ is a sub-algebra of $\mathrm{HS}(\operatorname{Inj_{coh}}(\mathbb{A}^{n+1},\mathbb{C}^{\ast},\omega))$. According to Proposition~\ref{extendedhochschild}, we have
 $$\Hom(\Delta,\Delta(t))\cong (\bigoplus_{g \in \mu_{d},\quad -\operatorname{rk}W_{g}\ \textbf{is even}}\Jac(\omega_{g})(t-k_{g}+d(\frac{-\operatorname{rk}W_{g}}{2})))^{\mathbb{C}^{\ast}}.$$
 In our case, if $g=1$, then $\operatorname{rk}W_{g}=0$, and $k_{g}=0$; if $g\neq 1$, then $\operatorname{rk}W_{g}=n+1$, $k_{g}=-n-1$, and $\mathrm{Jac}(\omega_{g})=k(0)$.
 Therefore, 
 \begin{itemize}
     \item If $n+1$ is odd, then 
     $\Hom(\Delta,\Delta(t))\cong  (\Jac(\omega)(t))^{\mathbb{C}^{\ast}}=\Jac(\omega)_{t}.$
     \item If $n+1$ is even, 
$$\Hom(\Delta,\Delta(t))\cong (\Jac(\omega)(t)\oplus E)^{\mathbb{C}^{\ast}}=\Jac(\omega)_{t}\oplus E^{\mathbb{C}^{\ast}},$$
     where $E=\bigoplus_{g\neq 1}\Jac(\omega_{g})(t+n+1-d(\frac{n+1}{2})) =\bigoplus_{g\neq 1}\Jac(\omega_{g})(t-\frac{(d-2)(n+1)}{2})$. If $t\neq \frac{(d-2)(n+1)}{2}$, then $E^{\mathbb{C}^{\ast}}=(\bigoplus^{d-1} k(t-\frac{(d-2)(n+1)}{2}))^{\mathbb{C}^{\ast}}=0$
      because $\Jac(\omega_{g})=k(0)$ for $g\neq 1$. Thus $\Hom(\Delta,\Delta(t))\cong \Jac(\omega)_{t}.$
 \end{itemize}

\end{proof}

\begin{rem}
In \cite{bondal2001reconstruction}, the canonical ring as the sub-algebra of Orlov's algebra $\mathrm{HA}$ described in section~\ref{section_orlov_algebra} is used to reconstruct smooth projective varieties with canonical line bundles ample or anti-ample. So it is reasonable to expect the subring of Serre algebra $\operatorname{Inj_{coh}}(\mathbb{A}^{n+1},\mathbb{C}^*,\omega)_S$ can be used to reconstruct the hypersurface defined by $\omega$. We will prove this expectation in Section~\ref{section_categorical_Torelli_hypersurfaces}. 
\end{rem}

\begin{rem}
The $\mathbb{Z}$-graded piece $\mathrm{Hom}(\Delta,\Delta(t))\cong\mathrm{HS}(\operatorname{Inj_{coh}}(\mathbb{A}^{n+1},\mathbb{C}^{\ast},\omega))(-k_1t,2k_2t+k_1(n+1)t)$.
\end{rem}

%\begin{thm}
%  If $gcd(2(n+1),d)=1$, then there is an algebra embedding
%  $$J(\omega)\subset Inj_{coh}(\mathbb{A}^{n+1},\mathbb{C}^{\ast},\omega)_{S}$$ 
%\end{thm}
%\begin{proof}
% We write $A$ for $Inj_{coh}(\mathbb{A}^{n+1},\mathbb{C}^{\ast},\omega)_{S}$. Take $i=2i'$, and $j=2j'$, then $(n+1)i+j$ is even, hence we have sub-algebra
% $$\bigoplus_{i',j'}Hom(\Delta,\Delta((d-2)(n+1)i'+dj').$$
% Since $gcd((d-2)(n+1),d)=gcd(2(n+1),d)=1$, we can choose integers $i'$, and $j'$ such that $1=(d-2)(n+1)i'+dj'$, and the same is true for any integer $t$.
% Then, by an appropriate choice of integers $i,j$, we have a sub-algebra
% $$\bigoplus^{(d-2)(n+1)}_{t=1}Hom(\Delta,\Delta(t))\subset A.$$
% According to \cite{}, $Hom(\Delta,\Delta(t))\cong Jac(\omega)_{t}$ if $t\geq 1$, and the isomorphism preserves the algebra structure.
%\end{proof}

\section{Categorical Torelli theorem for Fano hypersurfaces}\label{section_categorical_Torelli_hypersurfaces}
\begin{thm}\label{theorem_categorical_Torelli_thm_without_rotation_functors}
  Let $X$ and $X'$ be degree $d\leq n$ smooth hypersurfaces in $\mathbb{P}^n(n\geq 2)$ defined by $\omega$ and $\omega'$
  respectively. Assume $\operatorname{gcd}(n+1,d)=1$. If there is a Fourier-Mukai equivalence $\Ku(X)\simeq \Ku(X')$, then we have isomorphism $X\cong X'$.
\end{thm}

\begin{proof}
%We write $A$ for $Inj_{coh}(\mathbb{A}^{n+1},\mathbb{C}^{\ast},\omega)_{S}$. 

 First, if $n+1$ is odd, then the Serre algebra has a natural sub-algebra isomorphic to $\bigoplus_{t\geq 0} \Hom(\Delta, \Delta(t))$ by Proposition~\ref{subalgebra}. By the calculation in \cite[Theorem 5.39]{BFK}, the composition of $\bigoplus_{t\geq 0} \Hom(\Delta, \Delta(t))$ is the composition of functions under isomorphism $\bigoplus_{t\geq 0} \Hom(\Delta, \Delta(t))\cong \Jac(\omega)$, see also Remark~\ref{compatibility}. Therefore, $\Ku(X)\cong \Ku(X')$ implies an isomorphism of graded algebra $\Jac(\omega)\cong \Jac(\omega')$ by Theorem~\ref{Moritaserrealgebra}. Note that in degree one, it is a linear map $A$. Therefore, we have equality of ideal $\langle \partial_{i}(A\omega)\rangle=\langle \partial_{i}\omega'\rangle$. Then by Mather-Yau's reconstruction theorem\cite[Proposition 1.1]{Don}, $A\omega$ is projective equivalent to $\omega'$. Thus, $\omega$ is projective equivalent to $\omega'$, which implies $X\cong X'$. 
 \par
 Next, if $n+1$ is even, though we don't have $\Jac(\omega)$ as a natural sub-algebra, still we have natural graded piece $\Jac(\omega)_{d}$ since $d\neq \frac{(d-2)(n+1)}{2}$ by Proposition~\ref{subalgebra}, otherwise 
 $$d=\frac{(d-2)(n+1)}{2}>\frac{(d-2)d}{2}.$$
 Then $d\leq 3$. If $d=3$, then $n+1=6$, contradicts that $\operatorname{gcd}(n+1,d)=1$. It is clear that $d\neq 1$ and $d\neq 2$. Thus $\Ku(X)\cong \Ku(X')$ implies isomorphism $\Jac(\omega)_{d}\cong \Jac(\omega')_{d}$ which is induced by a linear transformation of degree one polynomials. Then similar argument as above implies $X\cong X'$.
\end{proof}

\begin{cor}
Categorical Torelli theorem holds for following Fano varieties:
\leavevmode\begin{enumerate}
    \item Cubic hypersurfaces of dimension $3k-1$ and $3k$ for $k\geq 1$. 
    \item Quintic fourfolds.
\end{enumerate}
\end{cor}

\begin{proof}
\leavevmode\begin{enumerate}
   \item Assume $d=3$, then $\mathrm{gcd}(3,n+1)=1$ implies $3\nmid n+1$, this means that dimension of $X$ is $n=3k-1$ or $3k$ with $k\geq 1$. 
    \item If $d=5$, and $n=5$, then $\mathrm{gcd}(6,5)=1$.
\end{enumerate}
Then the statement follows from Theorem~\ref{theorem_categorical_Torelli_thm_without_rotation_functors}.
\end{proof}

\section{Categorical Torelli Theorem with rotation functor}\label{categorical_torelli_with_rotation_functors}
In this section, we give a very simple proof of \cite[Theorem 1.2]{pirozhkov2022categorical} via matrix factorizations. 
\begin{thm}
  Let $X$ and $X'$ be smooth hypersurfaces of degree $d< n+1$ in $\mathbb{P}^{n}(n\geq 2)$. %Assume $(d,n+1)\neq (3,6)$. 
  If there is a Fourier-Mukai  equivalence of pairs $(\Ku(X),(1))\cong (\Ku(X'),(1)')$, then $X\cong X'$.   
\end{thm}
\begin{proof}
 Let $\omega$ and $\omega'$ define $X$ and $X'$ respectively. According to Theorem~\ref{prop_Orlov_identification}, there are  isomorphisms of pairs in $\operatorname{Hqe(dg-cat)}$
 \begin{align*}
  (\Ku_{dg}(X),(1))\cong &(\operatorname{Inj_{coh}}(\mathbb{A}^{n+1},\mathbb{C}^{\ast},\omega),\{1\})\\
   (\Ku_{dg}(X'),(1)')\cong &(\operatorname{Inj_{coh}}(\mathbb{A}^{n+1},\mathbb{C}^{\ast},\omega'),\{1\}')\\
 \end{align*}
 Which induces an isomorphism of pairs by diagram chasing
  $$(\operatorname{Inj_{coh}}(\mathbb{A}^{n+1},\mathbb{C}^{\ast},\omega),\{1\})\cong (\operatorname{Inj_{coh}}(\mathbb{A}^{n+1},\mathbb{C}^{\ast},\omega'),\{1\}').$$
  If $n+1$ is odd, then we have isomorphism of graded algebra by Proposition~\ref{subalgebra}and Theorem~\ref{Moritaserrealgebra},
  $$\bigoplus^{(n+1)(d-2)}_{t\geq 0}\Hom(\Delta,\Delta(t))\cong \bigoplus^{(n+1)(d-2)}_{t\geq 0} \Hom(\Delta',\Delta'(t)).$$
  That is, we have an  isomorphism of graded algebra
  $\Jac(\omega)\cong \Jac(\omega')$. Thus $\omega$ is projective equivalent to $\omega'$ by \cite[Proposition 1.1]{Don}.
  Note that while the statement in \cite[Proposition 1.1]{Don} involves the full Jacobian ring, the proof clearly uses only the $d$-th component.
  \par
  If $n+1$ is even, the case for $d=\frac{(d-2)(n+1)}{2}$ is $(d,n)=(3,5)$. But then $d-1\neq \frac{(d-2)(n+1)}{2}$. So in this case we have isomorphism 
  $\Jac(\omega)_{d-1}\cong \Jac(\omega')_{d-1}$ induced by linear map of degree one polynomials by Proposition~\ref{subalgebra} and Theorem~\ref{Moritaserrealgebra}, hence $\Jac(\omega)_{d}\cong \Jac(\omega')_{d}$ induced by linear map of degree one polynomials. 
  %Thus $X\cong X'$ by \cite[Proposition 1.1]{Don}. 
  Similarly, the case for $d-1=\frac{(d-2)(n+1)}{2}$ is $(d,n)=(3,3)$. So there is an isomorphism $\Jac(\omega)_{d}\cong \Jac(\omega')_{d}$ induced by linear map of degree one polynomials by Proposition~\ref{subalgebra} and Theorem~\ref{Moritaserrealgebra}, which again implies $X\cong X'$ by \cite[Proposition 1.1]{Don}.
\end{proof}

\begin{rem}\label{remark_two_ring_isomorphism}
Let $\{i\}:\operatorname{Inj_{coh}}(\mathbb{A}^{n+1},\mathbb{C}^{\ast},\omega)\simeq \operatorname{Inj_{coh}}(\mathbb{A}^{n+1},\mathbb{C}^{\ast},\omega)$ be the degree shift functor with corresponding Fourier-Mukai kernel $Q_i$. Define $L_{\mathrm{MF}}(X):=\bigoplus\mathrm{Hom}(Q_0,Q_i)$. On the other hand, the degree shift auto-equivalence $(i):\Ku(X)\simeq\Ku(X)$ is represented by Fourier-Mukai kernel $P_i$. Then we define another ring $L(X):=\bigoplus_i\mathrm{Hom}(P_0,P_i)$. In \cite{huybrechts2016hochschild}, the authors conjecture that $L_{\mathrm{MF}}(X)\cong L(X)$. Indeed, by \cite{pirozhkov2022categorical}, $\mathrm{Hom}(P_0,P_i)\cong\mathrm{Hom}(\mathrm{Id},(1)^i)$. Now since we have equivalence of the pair $\langle\Ku(X),(1)\rangle\simeq^{\phi} \langle \operatorname{Inj_{coh}}(\mathbb{A}^{n+1},G,\omega),\{1\}\rangle$, namely $(1)\cong \phi^{-1}\circ \{1\}\circ\phi$ by Theorem~\ref{prop_Orlov_identification}. We get
  $$\Hom(\Id,(1)^i)\cong \Hom(\Id,\{1\}^i).$$ Then 
$$L(X):=\bigoplus \mathrm{Hom}(P_{0},P_{i})\cong\bigoplus \mathrm{Hom}(Q_{0},Q_{i})\cong L_{MF}(X).$$ 
\end{rem}

%\section{Further comments}

\section{Appendix: Categorical Torelli theorem for weighted hypersurfaces}\label{section_appendix}
In this section, we illustrate the method used in proof of Theorem~\ref{theorem_categorical_Torelli_thm_without_rotation_functors} for a degree $6$ hypersurface in weighted projective space $\mathbb{P}(1,1,1,1,3)$, which is isomorphic to an index one prime Fano threefold of genus $2$. It is constructed as a double cover of $\mathbb{P}^3$ with branch divisor a sextic hypersurface. \emph{Categorical Torelli theorem} for this case was already established in \cite{Lahoz2023categorical} and \cite{Dell2023categorical} via completely different methods. We give a new proof.

%While we are preparing the paper, we are informed by Paolo Stellari that in \cite{Lahoz2023categorical}, they prove categorical Torelli theorem for degree $6$ hypersurfaces in weighted projective space $\mathbb{P}(1,1,1,1,3)$
%via a completely different method and the same statement is also proved in another upcoming preprint \cite{Dell2023categorical}. 

\begin{thm}\label{categorical_Torelli_weighted_hypersurfaces}
Let $X$ and $X'$ be smooth sextic hypersurfaces in weighted projective space $\mathbb{P}(1,1,1,1,3)$. Assume there is a Fourier-Mukai equivalence $\Ku(X)\simeq\Ku(X')$, then $X\cong X'$. 
\end{thm}

\begin{proof}
Consider Matrix Factorization $\operatorname{Inj_{coh}}(\mathbb{A}^{5},\mathbb{C}^{\ast},\omega)$, the weight of $\mathbb{C}^{\ast}$-action is $(1,1,1,1,3)$. According to \cite[Theorem 6.13]{BFK}, we have $\Ku(X)\cong \operatorname{Inj_{coh}}(\mathbb{A}^{5},\mathbb{C}^{\ast},\omega)$ and $\Ku(X')\cong Inj_{coh}(\mathbb{A}^{5},\mathbb{C}^{\ast},\omega')$, where $\omega$ and $\omega'$ are degree $6$ polynomial defining $X$ and $X'$ respectively. Then a Fourier-Mukai equivalence $\Ku(X)\cong \Ku(X')$ induces an equivalence $\operatorname{Inj_{coh}}(\mathbb{A}^{5},\mathbb{C}^{\ast},\omega)\cong \operatorname{Inj_{coh}(\mathbb{A}^{5},\mathbb{C}^{\ast},\omega)}$ in $\operatorname{Hqe(dg-cat)}$. Since $\operatorname{gcd}(\sum^{n+1}_{j}q_{j},d)=\operatorname{gcd}(7,6)=1$, according to Proposition~\ref{Moritaserrealgebra} and the same proof in Proposition~\ref{subalgebra}, we have isomorphism of algebra,
\begin{equation}\label{equation1}
\bigoplus_{t\geq 0}\Hom(\Delta,\Delta(t))\cong\bigoplus_{t\geq 0} \Hom(\Delta',\Delta'(t)).
\end{equation}
Then by \cite[Theorem 1.2] {BFK}, 
$$\Hom(\Delta,\Delta(t))\cong (\bigoplus_{g\in \mu_{6},\ -\operatorname{rk}W_{g}\  \textbf{is even}}\Jac(\omega_{g})(t-k_{g}+6(\frac{-\operatorname{rk}W_{g}}{2}))))^{\mathbb{C}^{\ast}}.$$
Write $\mu_{6}=\langle \lambda\rangle$. Then
$$(\mathbb{A}^{5})^{\lambda^{i}}=\begin{cases}
   (0,0,0,0,0);k_{\lambda^{i}}=-7; \operatorname{rk}(W_{\lambda^{i}})= 5,& if\  i=1, 3, 5,   \\
  (0,0,0,0,x_{5});k_{\lambda^{i}}=-4;\operatorname{rk}(W_{\lambda^{i}})=4, & if\  i=2, 4   \\  
  \mathbb{A}^{5};k_{\lambda^{i}}=0;\operatorname{rk}W_{\lambda^{i}}=0, & if\ i=6
\end{cases}$$
Write $\omega=x^{2}_{5}+f(x_{1},x_{2},x_{3},x_{4})$. Then $\Jac(\omega_{\lambda^{4}})=\Jac(\omega_{\lambda^{2}})=k[x_{5}]/\partial x^{2}_{5}=k(0)$, and $\Jac(\omega_{\lambda^{1}})=\Jac(\omega_{\lambda^{3}})=\Jac(\omega_{\lambda^{5}})=\Jac(\omega_{\lambda^{6}})=k(0)$. Therefore,
\begin{align*}
    \Hom(\Delta,\Delta(t))\cong\  &\bigoplus_{i}(\Jac(\omega_{\lambda^{i}})(t-k_{\lambda^{i}}+6(\frac{-\operatorname{rk}W_{\lambda^{i}}}{2})))^{\mathbb{C}^{\ast}}\\
    =\  & \Jac(\omega)_{t}\oplus k(t-8)\oplus k(t-8).
\end{align*}
Thus $\Hom(\Delta,\Delta(t))\cong \Jac(\omega)_{t}$ for $t\neq 8$. The same for $\Hom(\Delta',\Delta'(t))$, $\omega'=x^{2}_{5}+f'(x_{1},x_{2},{x_{3},x_{4}})$. According to isomorphism (1), we have $k[x_{1},x_{2},x_{3},x_{4}]_{6}/\langle\frac{\partial f}{\partial x_{i}}\rangle_{k}=\Jac(\omega)_{6}\cong \Jac(\omega')_{6}=k[x_{1},x_{2},x_{3},x_{4}]_{6}/\langle\frac{\partial f'}{\partial x_{i}}\rangle_{k}$  induced by autormorphism of degree one polynomials, which implies $f$ is projective equivalent to $f'$. Thus $X\cong X'$ by \cite[Proposition (1.3)]{Donagi1986}
\end{proof}

%\subsection{Statement of conflict of interest}
%On behalf of all authors, the corresponding author states that there is no conflict of interest.

% Of course, to make a completely rigorous argument, we should adopt the language of bi-modules of dg-categories as in \cite{polishchuk2014lefschetz}. We will do it later. 

%\bibliographystyle{spmpsci}
%\bibliography{Categoricaltorellihypersurface(1)}

%\bibliographystyle{alpha}
%{\small{\bibliography{Categoricaltorellihypersurface}}}

\end{document}